\documentclass{article}

\usepackage{amssymb}	
\usepackage{amsmath}    
\usepackage{setspace} 
\usepackage{latexsym}	
\usepackage{longtable}
\usepackage{multirow}
\usepackage{epsfig}
\allowdisplaybreaks

\newtheorem{algorithm}{Algorithm}[section]

\newtheorem{theorem}{Theorem}[section]

\newtheorem{remark}{Remark}[section]

\newcommand\bigzero{\makebox(0,0){\text{\huge0}}}

\def\A{{\bf A}}
\def\B{{\bf B}}

\def\E{{\bf E}}
\def\F{{\bf F}}

\def\H{{\bf H}}
\def\I{{\bf I}}

\def\0{{\bf 0}}
\def\P{{\bf P}}
\def\Q{{\bf Q}}
\def\R{{\bf R}}
\def\S{{\bf S}}
\def\T{{\bf T}}
\def\U{{\bf U}}

\def\W{{\bf W}}

\def\Y{{\bf Y}}
\def\Z{{\bf Z}}

\def\u{{\bf u}}

\def\x{{\bf x}}
\def\y{{\bf y}}
\def\z{{\bf z}}

\def\Tr{{\rm T}}

\def\diag{{\rm diag}}

\newcommand{\qed}{\nobreak \ifvmode \relax \else \ifdim\lastskip<1.5em \hskip-\lastskip \hskip1.5em plus0em minus0.5em \fi \nobreak \vrule height0.75em width0.5em depth0.25em\fi}

\title{An Efficient LQR Design for Discrete-Time Linear Periodic System Based on a Novel Lifting Method}

\author{Yaguang Yang\thanks{Instrumentation and Control Engineer, Office of Research, US NRC, Two White Flint North
11545 Rockville Pike, Rockville, MD 20852-2738. Phone: (301) 415-0655. Email: yaguang.yang@verizon.net.} }
 
\begin{document}

\maketitle

\begin{abstract}
This paper proposes a novel lifting method which converts the 
standard discrete-time linear periodic system to an augmented 
linear time-invariant system. The linear quadratic optimal 
control is then based on the solution of the discrete-time 
algebraic Riccati equation associated with the augmented linear 
time-invariant model. An efficient algorithm for solving the 
Riccati equation is derived by using the special structure of the 
augmented linear time-invariant system. It is shown that the 
proposed method is very efficient compared to the ones that use 
algorithms for discrete-time periodic algebraic Riccati equation. 
The efficiency and effectiveness of the proposed algorithm 
is demonstrated by the simulation test for the design problem of 
spacecraft attitude control using magnetic torques.
\end{abstract}

\newpage

\section{Introduction}

Many engineering systems are naturally periodic, for example, 
spacecraft attitude control using magnetic torques \cite{la04}, 
helicopter rotors control system \cite{abl00}, wind turbine
control system \cite{stol03}, networked control system 
\cite{zh06}, and multirate sampled data system \cite{ks91}.
It has been known for about six decades that linear periodic 
time-varying system can be converted to some equivalent linear 
time-invariant systems \cite{jm58,jm59}. The most popular and 
widely used methods that convert the linear periodic time-varying
model into linear time-invariant models are the so-called lifting
methods proposed in \cite{mb75,gl91}. Although these reduced 
linear time-invariant models are nice for analysis but they are
not very convenient for control system design. For example, 
the Linear Quadratic Regulator (LQR) design for linear periodic
system has been focused on the periodic system not on the 
equivalent linear time-invariant systems proposed in 
\cite{mb75,gl91}. This strategy leads to extensive research on 
the solutions of the periodic Riccati equations (see 
\cite{bcg86,bcn89,bittanti91,varga08,varga13} and references 
therein). For the discrete-time linear periodic system,
two efficient algorithms for Discrete-time Periodic Algebraic 
Riccati Equation (DPARE) are emerged \cite{hl94,yang17}.

In this paper, we propose a novel lifting method that converts
the linear periodic system to an augmented Linear Time-Invariant
(LTI) system. We show that the LQR design method can be directly 
applied to this LTI system. Moreover, by making full use of the
structure of the augmented LTI system, we can derive a very
efficient algorithm. We compare the new algorithm to the ones 
proposed in \cite{hl94,yang17}. In addition to some simple 
analysis on the efficiency, we demonstrate the efficiency and 
effectiveness of the new algorithm by the simulation test for 
the design problems of spacecraft attitude control using magnetic torques.

The remainder of the paper is organized as follows. Section 2 
briefly summarizes the algorithms of \cite{hl94,yang17} so that
we can compare the proposed algorithm to the existing ones and 
analyze the efficiency of these algorithms. Section 3 proposes 
a novel lifting method and applies some standard discrete-time
algebraic Riccati equation result to the augmented LTI model. 
This leads to a very efficient algorithm for the LQR design
for the linear periodic system. Section 4
demonstrates the efficiency and effectiveness of the algorithm 
by some numerical test. Conclusions are summarized in the last
section.

\section{Periodic LQR design based on linear periodic system}

In this section, we briefly review two efficient algorithms for
solving DPARE developed in \cite{hl94,yang17}. This will help us 
later in the comparison of the proposed method to the existing 
methods.

Let $p$ be an integer representing the total number of samples 
in one period in a periodic discrete-time system. We consider the 
following discrete-time linear periodic system given as follows:
\begin{eqnarray}
\x_{k+1}= \A_k \x_k +\B_k \u_k,
\label{discrete}
\end{eqnarray}
where $\A_k = \A_{k+p} \in \R^{n \times n}$ and 
$\B_k = \B_{k+p} \in \R^{n \times m}$ are periodic 
time-varying matrices. For this discrete-time linear periodic 
system (\ref{discrete}), the LQR state feedback control is to 
find the optimal $\u_k$ to minimize the following quadratic cost
function
\begin{equation}
\lim_{N \rightarrow \infty} \left( \min  \frac{1}{2} \x_N^{\Tr} \Q_N  \x_N
+\frac{1}{2} \sum_{k=0}^{N-1} \x_k^{\Tr} \Q_k  \x_k+ \u_k^{\Tr} \R_k  \u_k
\right)
\label{contCost}
\end{equation}
where 
\begin{eqnarray}
\Q_k= \Q_{k+p} \ge 0, \\
\R_k= \R_{k+p} > 0,
\end{eqnarray}
and the initial condition $\x_0$ is given. It is well-known that 
the LQR design for problem (\ref{discrete}-\ref{contCost}) can be
solved by using the periodic solution of the discrete-time 
periodic algebraic Riccati equation \cite{bittanti91}. Two 
efficient algorithms \cite{hl94,yang17} have been developed to 
solve $p$ $n$-dimensional matrix Riccati equations to find $p$
positive semidefinite matrices $\P_k$, $k=1, \ldots, p$. Given, 
$\P_k$, the periodic feedback controllers are given by the 
following equations:
\begin{equation}
\u_k = -(\R_k+\B_k^{\Tr}\P_{k}\B_k)^{-1} \B^{\Tr}_k  \P_{k} \A_k \x_k.
\label{optiSolu}
\end{equation}
We summarize these two algorithms as follows:
Let
\begin{equation}
\E_k =  \left[ \begin{array}{cc} 
\I & \B_k \R_k^{-1} \B_k^{\Tr}  \\ \0  & \A_{k}^{\Tr}
\end{array} \right]= \E_{k+p},
\label{Ek}
\end{equation}
\begin{equation}
\F_k =  \left[ \begin{array}{cc} 
\A_k & \0  \\ -\Q_k  & \I
\end{array} \right]= \F_{k+p}.
\label{Fk}
\end{equation}
If $\A_k$ is invertible, then $\E_k$ and $\F_k$ are invertible, 
and
\[
\E_k^{-1} =\left[ \begin{array}{cc}
\I & - \B_k \R_k^{-1} \B_k^{\Tr} \A_k^{-\Tr} \\
\0 & \A^{-\Tr}_k 
\end{array} \right]=\E_{k+p}^{-1}.
\]
and
\[
\F_k^{-1}=
\left[ \begin{array}{cc} 
\A_k^{-1}  & \0  \\ \Q_k  \A_k^{-1}   & \I
\end{array} \right]=\F_{k+p}^{-1}.
\]
Let $\y_k$ be the costate of $\x_k$, 
$\z_k =[ \x_k^{\Tr}, \y_k^{\Tr} ]^{\Tr}$, and 
\begin{equation}
\boldsymbol{\Pi}_k =\E_{k+p-1}^{-1} \F_{k+p-1} \E_{k+p-2}^{-1} \F_{k+p-2} \ldots \E_{k+1}^{-1} \F_{k+1} \E_k^{-1} \F_k =
\boldsymbol{\Pi}_{k+p},
\label{pi}
\end{equation}
\begin{equation}
\boldsymbol{\Gamma}_k = \F_{k}^{-1} \E_{k} \F_{k+1}^{-1} \E_{k+1} \ldots, 
\F_{k+p-2}^{-1} \E_{k+p-2} \F_{k+p-1}^{-1} \E_{k+p-1} =
\boldsymbol{\Gamma}_{k+p}.
\label{Gamma}
\end{equation}
The solutions of $p$ discrete-time periodic algebraic Riccati 
equations are symmetric positive semi-definite matrices, $\P_k$, 
$k=1, \ldots, p$, which are related to the solutions of either one
of the two linear systems of equations \cite{hl94,yang17}:
\begin{equation}
\z_{k+p}= \boldsymbol{\Pi}_k \z_k,
\label{lump}
\end{equation}
\begin{equation}
\z_k=\boldsymbol{\Gamma}_k \z_{k+p}.
\label{inverse}
\end{equation}
Therefore, $\P_k$, $k=1, \ldots, p$, can be obtained by two 
methods. The first method uses Schur decomposition:
\begin{equation}
\left[ \begin{array}{cc} \T_{11k} &  \T_{12k}  \\
\T_{21k}   &    \T_{22k}
\end{array} \right]^{\Tr}\boldsymbol{\Pi}_k
\left[ \begin{array}{cc} \T_{11k} &  \T_{12k}  \\
\T_{21k}   &    \T_{22k}
\end{array} \right]=
\left[ \begin{array}{cc}
\S_{11k} &  \S_{12k}  \\
\0   &    \S_{22k}
\end{array} \right],
\label{schur}
\end{equation}
where $\S_{11k}$ is upper-triangular and has all of its eigenvalues inside the unique circle. The periodic solution $\P_k$, $k=1,\ldots,p$, is given by \cite{hl94}
\begin{equation}
\P_{k}= \T_{21k}\T_{11k}^{-1}.
\label{Xk}
\end{equation}
The second method uses Schur decomposition:
\begin{equation}
\left[ \begin{array}{cc} \W_{11k} &  \W_{12k}  \\
\W_{21k}   &    \W_{22k}
\end{array} \right]^{\Tr}\boldsymbol{\Gamma}_k 
\left[ \begin{array}{cc} \W_{11k} &  \W_{12k}  \\
\W_{21k}   &    \W_{22k}
\end{array} \right]=
\left[ \begin{array}{cc}
\U_{11k} &  \U_{12k}  \\
\0   &    \U_{22k}
\end{array} \right],
\label{decomp1}
\end{equation}
where $\W_{11k}$ is upper-triangular and has all of its eigenvalues outside the unique circle. The periodic solution $\P_k$, $k=1,\ldots,p$, is given by \cite{yang17}
\begin{equation}
\P_k =\W_{21k}\W_{11k}^{-1}.
\label{solution}
\end{equation}
\begin{remark}
When $\A_k$ and $\Q_k$ are constant matrices, the second method is
much efficient because $\F_k$ becomes a constant matrix and
$\F_{k}^{-1} = \cdots = \F_{k+p-1}^{-1} = \F^{-1}$, which 
makes the computation of (\ref{Gamma}) much more efficient than the
computation of (\ref{pi}).
\end{remark}

\section{Periodic LQR design based on linear time-invariant system}

We propose a lifting method in this section to convert the 
discrete-time linear periodic system into an augmented linear 
time-invariant system. Thereby, the periodic LQR design is 
reduced to the LQR design for the augmented linear time-invariant 
system.

To simplify our discussion, let us consider a periodic system with
$p=3$. We will use $k$ for the discrete-time in the periodic 
system and $K$ for the discrete-time in the augmented system.
\begin{subequations}
\begin{gather}
\x_1=\A_0 \x_0 + \B_0 \u_0,  \nonumber \\
\x_2=\A_1 \x_1 + \B_1 \u_1,  \nonumber \\
\x_3=\A_2 \x_2 + \B_2 \u_2,  \nonumber \\
\x_4=\A_0 \x_3 + \B_0 \u_3,  \nonumber \\
\x_5=\A_1 \x_4 + \B_1 \u_4,  \nonumber \\
\x_6=\A_2 \x_5 + \B_2 \u_5,  \nonumber \\
\x_7=\A_0 \x_6 + \B_0 \u_6,  \nonumber \\
\vdots  \nonumber 
\end{gather}
\label{simpleExample}
\end{subequations}
We can easily regroup the periodic system and rewrite
it as the following form:
\begin{eqnarray}
\bar{\x}_1 & =& \left[ \begin{array}{c}
\x_1  \\  \x_2  \\  \x_3
\end{array} \right]
=\left[ \begin{array}{ccc}
\0  &  \0  &  \A_0   \\
\0  &  \0  &  \A_1 \A_0   \\
\0  &  \0  &   \A_2  \A_1 \A_0   
\end{array} \right]
\left[ \begin{array}{c}
\0  \\  \0  \\  \x_0
\end{array} \right]           +
\left[ \begin{array}{ccc}
\B_0  &  \0  &  \0   \\
\A_1 \B_0  &  \B_1  &  \0   \\
\A_2 \A_1 \B_0    &  \A_2 \B_1  &   \B_2   
\end{array} \right]
\left[ \begin{array}{c}
\u_0  \\  \u_1  \\  \u_2
\end{array} \right]  \nonumber \\
 & =& \bar{\A} \bar{\x}_0 + \bar{\B} \bar{\u}_0,  \nonumber
\label{simple1}
\end{eqnarray}
\begin{eqnarray}
\bar{\x}_2 & =&  \left[ \begin{array}{c}
\x_4  \\  \x_5  \\  \x_6
\end{array} \right]
=\left[ \begin{array}{ccc}
\0  &  \0  &  \A_0   \\
\0  &  \0  &  \A_1 \A_0   \\
\0  &  \0  &   \A_2  \A_1 \A_0   
\end{array} \right]
\left[ \begin{array}{c}
\x_1  \\  \x_2  \\  \x_3
\end{array} \right]           +
\left[ \begin{array}{ccc}
\B_0  &  \0  &  \0   \\
\A_1 \B_0  &  \B_1  &  \0   \\
\A_2 \A_1 \B_0    &  \A_2 \B_1  &   \B_2   
\end{array} \right]
\left[ \begin{array}{c}
\u_3  \\  \u_4  \\  \u_5
\end{array} \right]  \nonumber \\
 & =& \bar{\A} \bar{\x}_1 + \bar{\B} \bar{\u}_1,  \nonumber
\label{simple2}
\end{eqnarray}
in general, for $k \ge 0$ ($K \ge 0$), we have
\begin{eqnarray}
& & \bar{\x}_{K+1}   \nonumber \\
& := & \left[ \begin{array}{c}
\x_{pk+1}  \\  \x_{pk+2}  \\  \x_{pk+3}
\end{array} \right]
=\left[ \begin{array}{ccc}
\0  &  \0  &  \A_0   \\
\0  &  \0  &  \A_1 \A_0   \\
\0  &  \0  &   \A_2  \A_1 \A_0   
\end{array} \right]
\left[ \begin{array}{c}
\x_{p(k-1)+1}  \\  \x_{p(k-1)+2}  \\  \x_{p(k-1)+3}
\end{array} \right]           +
\left[ \begin{array}{ccc}
\B_0  &  \0  &  \0   \\
\A_1 \B_0  &  \B_1  &  \0   \\
\A_2 \A_1 \B_0    &  \A_2 \B_1  &   \B_2   
\end{array} \right]
\left[ \begin{array}{c}
\u_{pk}  \\  \u_{pk+1}  \\  \u_{pk+2}
\end{array} \right]  \nonumber \\
 & :=& \bar{\A} \bar{\x}_{K} + \bar{\B} \bar{\u}_{K},  
\label{simple3}
\end{eqnarray}
where 
\[
\bar{\x}_{0} = \left[ \begin{array}{c}
\x_{-2}  \\  \x_{-1}  \\  \x_{0}
\end{array} \right] 
:= \left[ \begin{array}{c}
\0  \\  \0  \\  \x_{0}
\end{array} \right] , \hspace{0.1in}
\bar{\u}_{0} = \left[ \begin{array}{c}
\u_{0}  \\  \u_{1}  \\  \u_{2}
\end{array} \right].
\]
It is worthwhile to note that (\ref{simple3}) is a linear 
time-invariant system.
We can easily extend the result to the general case. Let
\begin{equation}
\bar{\x}_{K}=\left[  \begin{array}{c}
\x_{p(k-1)+1} \\ \x_{p(k-1)+2} \\  \vdots  \\   \x_{p(k-1)+p} 
\end{array} \right], \hspace{0.1in} 
\bar{\x}_{0}:=\left[  \begin{array}{c}
\0   \\  \vdots  \\ \0 \\   \x_{0} 
\end{array} \right], \hspace{0.1in}
\bar{\u}_{K}= \left[  \begin{array}{c}
\u_{pk} \\  \u_{pk+1} \\  \ldots \\ \u_{pk+p-1}
\end{array} \right]. \nonumber
\end{equation}
We will use the following notation:
\begin{equation}
\bar{\x}_{K+1}=\left[  \begin{array}{c}
\x_{pk+1} \\ \x_{pk+2} \\  \vdots  \\   \x_{pk+p} 
\end{array} \right],     \hspace{0.1in} 
\bar{\u}_{K+1}= \left[  \begin{array}{c}
\u_{p(k+1)} \\  \u_{p(k+1)+1} \\  \ldots \\ \u_{p(k+1)+p-1}
\end{array} \right]. \nonumber
\end{equation}

\begin{theorem}
Given a linear periodic discrete-time system with period of $p$
as follows:
\begin{eqnarray}
\x_{pk+1} & = & \A_0 \x_{p(k-1)+1} + \B_0 \u_{pk} ,  \nonumber \\
\x_{pk+2} & = & \A_1 \x_{p(k-1)+2} + \B_1 \u_{pk+1},  \nonumber \\
& \vdots  & \nonumber \\
\x_{pk+p} & = & \A_{p-1} \x_{p(k-1)+p} + \B_{p-1}  \u_{pk+p-1}. 
\label{generalCase}
\end{eqnarray}
Then, this discrete-time periodic system is equivalent to the 
linear time-invariant system given as follows:
\begin{eqnarray}
\bar{\x}_{K+1} :=
\left[ \begin{array}{c}
\x_{pk+1}  \\  \x_{pk+2}  \\  \vdots  \\  \x_{pk+p}
\end{array} \right]
& = & \left[ \begin{array}{cccc}
\0  &  \ldots &    \0  &            \A_0   \\
\0  &  \ldots &    \0  &    \A_1 \A_0   \\
\vdots &  \vdots &  \vdots &  \vdots   \\
\0  &  \ldots &    \0  &   \A_{p-1} \ldots  \A_2 \A_1 \A_0   
\end{array} \right]
\left[ \begin{array}{c}
\x_{p(k-1)+1}  \\  \x_{p(k-1)+2}  \\ \vdots \\  \x_{p(k-1)+p}
\end{array} \right]          \nonumber \\
& + &
\left[ \begin{array}{cccc}
\B_0  &  \0   &  \ldots     &  \0   \\
\A_1 \B_0  &  \B_1  &  \ldots     &  \0   \\
\vdots  &  \vdots  & \vdots   &   \vdots  \\
 \A_{p-1} \ldots  \A_1 \B_0    & \A_{p-1} \ldots \A_2 \B_1  &     \ldots &    \B_{p-1}   
\end{array} \right]
\left[ \begin{array}{c}
\u_{pk}  \\  \u_{pk+1} \\   \ldots  \\  \u_{pk+p-1}
\end{array} \right]          \nonumber \\
& := & \bar{\A} \bar{\x}_{K} + \bar{\B} \bar{\u}_{K},
\label{linearTimeInvariant}
\end{eqnarray}
where $\bar{\A} \in \R^{pn \times pn}$ and 
$\bar{\B} \in \R^{pn \times pm}$. Moreover, the structure of 
$\bar{\B}$ matrix guarantees the causality of the system 
(\ref{linearTimeInvariant}).
\hfill \qed
\label{liftingTheorem}
\end{theorem}

It is worthwhile to emphasize that there is no overlap
between $\bar{\x}_{K+1}$ and $\bar{\x}_{K}$; in addition,
there is no overlap between $\bar{\u}_{K+1}$ and 
$\bar{\u}_{K}$. This is the major difference between the proposed 
lifting method and the existing lifting methods in 
\cite{mb75,gl91} (see also \cite{varga13}). This feature makes it 
possible to apply existing design methods for linear 
time-invariant system (\ref{linearTimeInvariant}) which is 
equivalent to the linear periodic system (\ref{generalCase}). 
In the remainder of the paper, we
will investigate the LQR design for the system 
(\ref{linearTimeInvariant}). The LQR state feedback control is to 
find the optimal $\bar{\u}_{K}$ to minimize the following quadratic cost
function
\begin{equation}
\lim_{N \rightarrow \infty} \left( \min  \frac{1}{2} 
\bar{\x}_N^{\Tr} \bar{\Q}_N  \bar{\x}_N
+\frac{1}{2} \sum_{K=0}^{N-1} \bar{\x}_K^{\Tr} 
\bar{\Q}_K  \bar{\x}_K+ \bar{\u}_{K}^{\Tr} 
\bar{\R}_K  \bar{\u}_{K}
\right)
\label{augmentedCost}
\end{equation}
where 
\begin{eqnarray}
\bar{\Q}_K := \diag ({\Q}_0, \ldots, \Q_{p-1}) \ge 0, \hspace{0.1in} 
\bar{\R}_K := \diag ({\R}_0, \ldots, \R_{p-1}) > 0,
\label{QRblock}
\end{eqnarray}
are constant matrices and the initial condition $\bar{\x}_0$ is 
given. It is straightforward to see that the optimal control problem 
described by (\ref{linearTimeInvariant}) and (\ref{augmentedCost})
is time-invariant but equivalent to the time-varying periodic system
described by (\ref{discrete}) and (\ref{contCost}). Moreover, the 
optimal feedback matrix of the system (\ref{linearTimeInvariant}
-\ref{augmentedCost}) is given as follows:
\begin{equation}
\bar{\u}_K= - (\bar{\R}+\bar{\B}^{\Tr} \bar{\P} \bar{\B})^{-1}
\bar{\B}^{\Tr} \bar{\P} \bar{\A} \bar{\x}_K,
\label{augmentedU}
\end{equation}
where $\bar{\P}$ is the solution of the following {\it time-invariant 
algebraic Riccati equation}:
\begin{equation}
\bar{\A}^{\Tr} \bar{\P}\bar{\A} -\bar{\P} -\bar{\A}^{\Tr} \bar{\P}
\bar{\B} (\bar{\R}+\bar{\B}^{\Tr} \bar{\P} \bar{\B})^{-1}
\bar{\B}^{\Tr} \bar{\P} \bar{\A} + \bar{\Q} = \0.
\label{TIriccati}
\end{equation} 
Notice that $\bar{\A}$ is not invertible, this algebraic Riccati 
equation cannot be directly solved by using the algorithms for 
time-invariant algebraic Riccati equation proposed in 
\cite{vaughan70,laub79}, but it can be solved using the algorithm 
proposed in \cite{pls80}. In this paper, we propose a more 
efficient algorithm than the one of \cite{pls80}. The new
algorithm makes full use of the specific structure 
of $\bar{\A}$ in which the first $(p-1)n$ columns are zeros. 
Denote 
\begin{equation}
\bar{\Q}:=\bar{\Q}_K
=\left[ \begin{array}{c|c}
\diag( {\Q}_0, \ldots, {\Q}_{p-2}) & \0  \\
\hline
\0  &  {\Q}_{p-1}   \end{array}   \right]
=\diag( \bar{\Q}_1, \bar{\Q}_{2}  ),
\label{Qblock}
\end{equation}
where $\bar{\Q}_1=\diag( {\Q}_0, \ldots, {\Q}_{p-2})  
\in \R^{(p-1)n \times (p-1)n}$ 
and $\bar{\Q}_{2}= {\Q}_{p-1} \in \R^{n \times n}$,
\begin{equation}
\bar{\R}:=\bar{\R}_K
=\left[ \begin{array}{c|c}
\diag( {\R}_0, \ldots, {\R}_{p-2}) & \0  \\
\hline
\0  &  {\R}_{p-1}   \end{array}   \right]
=\diag( \bar{\R}_1, \bar{\R}_{2}  ),
\label{Rblock}
\end{equation}
where $\bar{\R}_1=\diag( {\R}_0, \ldots, {\R}_{p-2}) 
\in \R^{(p-1)m \times (p-1)m}$ 
and $\bar{\R}_{2}= {\R}_{p-1} \in \R^{m \times m}$.  Let
\begin{eqnarray}
\bar{\A} &  =  &
\left[ \begin{array}{ccc|c}
\0  &  \ldots &    \0  &            \A_0   \\
\0  &  \ldots &    \0  &    \A_1 \A_0   \\
\vdots &  \vdots &  \vdots &  \vdots   \\
\0  &  \ldots &    \0  &   \A_{p-2} \ldots  \A_1 \A_0   \\
\hline
\0  &  \ldots &    \0  &   \A_{p-1} \ldots  \A_2 \A_1 \A_0   
\end{array} \right]
=\left[ 
\begin{smallmatrix}
\0  &  \ldots &    \0    \\
\0  &  \ldots &    \0   \\
\vdots &  \vdots &  \vdots  \\
\0  &  \ldots &    \0   \\
\0  &  \ldots &    \0     
\end{smallmatrix} \middle|
\begin{smallmatrix}
  \A_0   \\
\A_1 \A_0   \\
\vdots   \\
\A_{p-2} \ldots \A_1 \A_0   \\
\hline
\A_{p-1} \ldots  \A_2 \A_1 \A_0 
\end{smallmatrix}
\right]   \nonumber \\
& = & \left[ \begin{array}{c|c}
\underbrace{
\begin{smallmatrix}
\\ \\
\bigzero
\\ \\  \\  \\
\end{smallmatrix}
}_{(p-1)n \mbox{   columns}}  
& 
\underbrace{
\begin{smallmatrix}
\bar{\A}_{1} \\
\\
\hline
\\ \\  
\bar{\A}_{2}
\\ \\
\end{smallmatrix} 
}_{n \mbox{   columns}} 
\end{array} \right]
=\left[
\begin{smallmatrix} \0  \end{smallmatrix} \middle| 
\begin{smallmatrix}  \bar{\F}  \end{smallmatrix} 
\right],
\label{Ablock}
\end{eqnarray}
where $\bar{\A}_{1} \in \R^{(p-1)n \times n}$
$\bar{\A}_{2} \in \R^{n \times n}$, 
$\bar{\F} =\left[ \bar{\A}_1^{\Tr}, \bar{\A}_2^{\Tr} \right]^{\Tr} 
\in \R^{pn \times n}$,
\begin{eqnarray}
\bar{\B} = \left[ \begin{array}{ccccc}
\B_0  &  \0   &  \ldots     &  \0      &  \0  \\
\A_1 \B_0  &  \B_1  &  \ldots     &  \0      &  \0  \\
\vdots  &  \vdots  & \vdots   &   \vdots    &   \vdots  \\
\A_{p-2} \ldots  \A_1 \B_0    & \A_{p-2} \ldots \A_2 \B_1  &     \ldots &    \B_{p-2}      &  \0 \\
\hline
 \A_{p-1} \ldots  \A_1 \B_0    & \A_{p-1} \ldots \A_2 \B_1  &     \ldots &    \A_{p-1} \B_{p-2} & \B_{p-1}   
\end{array} \right]
= \left[ \begin{array}{c}
\\ 
\bar{\B}_1   \\
\\ 
\hline
\\
 \bar{\B}_{2}   
\end{array} \right],
\label{Bblock}
\end{eqnarray}
where $\bar{\B}_1 \in \R^{(p-1)n \times pm}$ and 
$\bar{\B}_{2} \in \R^{n \times pm}$, and
\begin{equation}
\bar{\P} = \left[ \begin{array}{cc}
\bar{\P}_{11} & \bar{\P}_{12}  \\
\bar{\P}_{21} & \bar{\P}_{22}  
\end{array} \right],
\label{Pblock}
\end{equation}
where $\bar{\P}_{11} \in \R^{(p-1)n \times (p-1)n}$, 
$\bar{\P}_{12} \in \R^{(p-1)n \times n}$, 
$\bar{\P}_{21} \in \R^{n \times (p-1)n}$, and 
$\bar{\P}_{22} \in \R^{n \times n}$.
Let 
\begin{equation}
\Y = \bar{\P}
\bar{\B} (\bar{\R}+\bar{\B}^{\Tr} \bar{\P} \bar{\B})^{-1}
\bar{\B}^{\Tr} \bar{\P}.
\label{Yblock}
\end{equation}
Substituting (\ref{Qblock}), (\ref{Rblock}), (\ref{Ablock}), (\ref{Bblock}),
(\ref{Pblock}), and (\ref{Yblock}) into (\ref{TIriccati}) yields
\begin{equation}
\left[ \begin{array}{c}
\0  \\  \bar{\F}^{\Tr}  
\end{array} \right]
\bar{\P}  \left[ \begin{array}{cc}
\0  &  \bar{\F}  
\end{array} \right] 
-\left[ \begin{array}{cc}
\bar{\P}_{11} & \bar{\P}_{12}  \\
\bar{\P}_{21} & \bar{\P}_{22}  
\end{array} \right]
-\left[ \begin{array}{c}
\0  \\  \bar{\F}^{\Tr}  
\end{array} \right]
\bar{\Y}  \left[ \begin{array}{cc}
\0  &  \bar{\F}  
\end{array} \right] 
+ \left[ \begin{array}{cc}
\bar{\Q}_{1} & \0  \\
\0 & \bar{\Q}_{2}  
\end{array} \right] = \0,
\end{equation}
or equivalently
\begin{equation}
\left[ \begin{array}{cc}
\0  &  \0   \\
\0  &  \bar{\F}^{\Tr}  \bar{\P}   \bar{\F}
\end{array} \right]
-\left[ \begin{array}{cc}
\bar{\P}_{11} & \bar{\P}_{12}  \\
\bar{\P}_{21} & \bar{\P}_{22}  
\end{array} \right]
-\left[ \begin{array}{cc}
\0  &  \0   \\
\0  &  \bar{\F}^{\Tr}  \bar{\Y}   \bar{\F}
\end{array} \right] 
+ \left[ \begin{array}{cc}
\bar{\Q}_{1} & \0  \\
\0 & \bar{\Q}_{2}  
\end{array} \right] = \0.
\label{riccatiBlock}
\end{equation}
This proves $\bar{\P}_{12}=\bar{\P}_{21}^{\Tr}=\0$
and $\bar{\P}_{11} = \bar{\P}_{11}^{\Tr} = \bar{\Q}_{1}$. 
By examining the lower right block of (\ref{riccatiBlock}),
we have
\begin{equation}
\bar{\F}^{\Tr} \bar{\P} \bar{\F} = \bar{\A}_{1}^{\Tr} \bar{\Q}_1 
\bar{\A}_1 + \bar{\A}_{2}^{\Tr} \bar{\P}_{22} \bar{\A}_2
\in \R^{n \times n},
\label{lastBlock1}
\end{equation}
and
\begin{eqnarray}
& & \bar{\F}^{\Tr} \bar{\Y} \bar{\F} \nonumber \\
& = & \left[  \begin{array}{cc}
\bar{\A}_{1}^{\Tr}, \bar{\A}_{2}^{\Tr} 
\end{array} \right]
\left[  \begin{array}{c}
\bar{\Q}_1 \bar{\B}_{1} \\  \bar{\P}_{22} \bar{\B}_2
\end{array} \right]
\left[  \begin{array}{cc}
\bar{\R} + \bar{\B}_{1}^{\Tr} \bar{\Q}_1 \bar{\B}_{1}
+ \bar{\B}_{2}^{\Tr}  \bar{\P}_{22}^{\Tr} \bar{\B}_{2}
\end{array} \right]^{-1}
\left[  \begin{array}{cc}
\bar{\B}_{1}^{\Tr} \bar{\Q}_1 & 
\bar{\B}_2^{\Tr} \bar{\P}_{22}
\end{array} \right]
\left[  \begin{array}{c}
\bar{\A}_{1} \\ \bar{\A}_{2}
\end{array} \right]  \nonumber \\
& = & \underbrace{ \left[  
\bar{\A}_{1}^{\Tr} \bar{\Q}_1 \bar{\B}_{1} 
+ \bar{\A}_{2}^{\Tr} \bar{\P}_{22} \bar{\B}_2
\right] 
}_{n \times pm}
\underbrace{ \left[  \begin{array}{cc}
\bar{\R} + \bar{\B}_{1}^{\Tr} \bar{\Q}_1 \bar{\B}_{1}
+ \bar{\B}_{2}^{\Tr}  \bar{\P}_{22} \bar{\B}_{2}
\end{array} \right]^{-1}
}_{pm \times pm}
\underbrace{ \left[  
\bar{\B}_{1}^{\Tr} \bar{\Q}_1 \bar{\A}_{1} 
+ \bar{\B}_2^{\Tr} \bar{\P}_{22} \bar{\A}_{2}
\right]
}_{pm \times n}.
\label{lastBlock2}
\end{eqnarray}
Let 
\begin{subequations}
\begin{gather}
\hat{\A} = \bar{\A}_{2} \in \R^{n \times n}, \\
\hat{\B} = \bar{\B}_2 \in \R^{n \times pm}, \\
\hat{\Q} = \bar{\Q}_2 + \bar{\A}_{1}^{\Tr} \bar{\Q}_1 \bar{\A}_1
 \in \R^{n \times n},  \\
\hat{\R} = \bar{\R} + \bar{\B}_{1}^{\Tr} \bar{\Q}_1 \bar{\B}_1
 \in \R^{pm \times pm}, \\
\hat{\S} = \bar{\A}_{1}^{\Tr} \bar{\Q}_1 \bar{\B}_{1}
 \in \R^{n \times pm},  \\
\hat{\P} = \bar{\P}_{22}  \in \R^{n \times n}.
\end{gather}
\label{let1}
\end{subequations}
We can rewrite the lower right block of (\ref{riccatiBlock}) as 
follows:
\begin{equation}
\hat{\A}^{\Tr} \hat{\P} \hat{\A} - \hat{\P} 
- \underbrace{ \left( \hat{\A}^{\Tr} \hat{\P} \hat{\B} + \hat{\S}
 \right) }_{n \times pm}
\underbrace{ \left( \hat{\B}^{\Tr} \hat{\P} \hat{\B} +\hat{\R} 
\right)^{-1} }_{pm \times pm}
\underbrace{ \left( \hat{\B}^{\Tr} \hat{\P} \hat{\A} 
+ \hat{\S}^{\Tr} \right) }_{n \times pm}
+ \hat{\Q} = \0.
\label{reducedRiccati}
\end{equation}
The Riccati equation (\ref{reducedRiccati}) is a special case
discussed in \cite{al84}. An efficient Matlab function {\tt dare} 
that implements an algorithm of \cite{al84} is available to solve
(\ref{reducedRiccati}).

\begin{remark}
Comparing to the methods described in the previous section which need to solve $p$ $n$-dimensional discrete-time Riccati 
equations, we need only to solve one $n$-dimensional 
discrete-time Riccati equation using the method proposed in this 
section. 
\end{remark}

To compare the efficiency of the method to the ones developed in 
\cite{hl94,yang17}, we will not use the Matlab function {\tt dare}
because {\tt dare} calculates more information than the solution of the
Riccati equation (\ref{reducedRiccati}). Let $\tilde{\B} =\hat{\B}$,
$\tilde{\R} =\hat{\R}$, 
\begin{equation}
\tilde{\A} = \hat{\A} -\hat{\B} \hat{\R}^{-1}\hat{\S}^{\Tr},
\label{tiltA}
\end{equation}
and
\begin{equation}
\tilde{\Q} = \hat{\Q} -\hat{\S} \hat{\R}^{-1}\hat{\S}^{\Tr}.
\label{tiltQ}
\end{equation}
Riccati equation (\ref{reducedRiccati}) can be solved by either
eigen-decomposition or Schur decomposition for the following
generalized eigenvalue problem \cite[page 1748, equation (8)]{al84}:
\begin{equation}
{\lambda} \left[ \begin{array}{cc}
\I & \tilde{\B} \tilde{\R}^{-1}\tilde{\B}^{\Tr}\\
\0 & \tilde{\A} ^{\Tr}   \end{array} \right]
- \left[ \begin{array}{cc}
\tilde{\A} & \0 \\
-\tilde{\Q} &  \I   \end{array} \right]
:= \lambda \E - \F.
\label{gEigen}
\end{equation}
If $\tilde{\A}$ is invertible, the problem is reduced to solve
the following eigenvalue problem \cite[equation (30)]{yang17}:
\begin{eqnarray}
\Z = \F^{-1} \E = 
 \left[ \begin{array}{cc} 
\tilde{\A}^{-1}  & \tilde{\A}^{-1} \tilde{\B}\tilde{\R}^{-1}\tilde{\B}^{\Tr}  
\\ \tilde{\Q} \tilde{\A}^{-1}   &  \tilde{\A}^{\Tr}+\tilde{\Q} \tilde{\A}^{-1} 
\tilde{\B}\tilde{\R}^{-1}\tilde{\B}^{\Tr}
\end{array} \right].
\label{reducedMatrix}
\end{eqnarray}
Using Schur decomposition for (\ref{reducedMatrix}), we have
\begin{equation}
\left[ \begin{array}{cc} \W_{11} &  \W_{12}  \\
\W_{21}   &    \W_{22}
\end{array} \right]^{\Tr}\boldsymbol{\Z} 
\left[ \begin{array}{cc} \W_{11} &  \W_{12}  \\
\W_{21}   &    \W_{22}
\end{array} \right]=
\left[ \begin{array}{cc}
\S_{11} &  \S_{12}  \\
\0   &    \S_{22}
\end{array} \right],
\label{decompZ}
\end{equation}
where $\S_{11}$ is upper-triangular and has all of its eigenvalues 
outside the unique circle. The solution of the discrete algebraic Riccati
equation (\ref{reducedRiccati}) is given by
\begin{equation}
\hat{\P} =\W_{21}\W_{11}^{-1}.
\label{Psolution}
\end{equation}

We summarize the proposed algorithm as follows:

\begin{algorithm} {\ } \\
Data: $\A_0, \ldots, \A_{p-1}$, $\B_0, \ldots, \B_{p-1}$,
$\Q_0, \ldots, \Q_{p-1}$, $\R_0, \ldots, \R_{p-1}$.
\begin{itemize}
\item[] Step 1: Form 
\begin{subequations}
\begin{gather}
\bar{\A}_1 = \left[ \begin{array}{c}
\A_0   \\
\A_1 \A_0   \\
\vdots   \\
\A_{p-2} \ldots  \A_2 \A_1 \A_0   
\end{array} \right], \hspace{0.1in}
\bar{\B_1} = \left[ \begin{array}{ccccc}
\B_0  &  \0   &  \ldots     &  \0 &  \0     \\
\A_1 \B_0  &  \B_1  &  \ldots     &  \0  &  \0    \\
\vdots  &  \vdots  & \vdots   &   \vdots  &   \vdots   \\
 \A_{p-2} \ldots  \A_1 \B_0,    & \A_{p-2} \ldots \A_2 \B_1,  &     \ldots &    \B_{p-2},   &  \0   
\end{array} \right],
\label{AB1}
\\
\bar{\A}_2 = \A_{p-1} \ldots  \A_2 \A_1 \A_0, \hspace{0.1in}
\bar{\B}_2 =  \left[ \begin{array}{cccc}
\A_{p-1} \ldots  \A_1 \B_0,    & \A_{p-1} \ldots 
\A_2 \B_1,  &     \ldots &    \B_{p-1}   
\end{array} \right],
\\
\bar{\Q}_1=\diag( {\Q}_0, \ldots, {\Q}_{p-2}), \hspace{0.1in}
\bar{\Q}_{2}= {\Q}_{p-1},
\\
\bar{\R}_1=\diag( {\R}_0, \ldots, {\R}_{p-2}), \hspace{0.1in}
\bar{\R}_{2}= {\R}_{p-1}.
\end{gather}
\label{AB2}
\end{subequations}
\item[] Step 2: Form $\hat{\A}$, $\hat{\B}$, $\hat{\Q}$, 
$\hat{\R}$, and $\hat{\S}$ using (\ref{let1}).
\item[] Step 3: Find the solution $\hat{\P}$ of the discrete-time 
algebraic Riccati equation (\ref{reducedRiccati}) using the 
algorithm of \cite{al84} implemented as {\tt dare} or using the 
algorithm of \cite{yang17}, i.e., (\ref{decompZ}) and 
(\ref{Psolution}).
\item[] Step 4: The solution of the discrete-time algebraic 
Riccati equation (\ref{TIriccati}) is given by 
\begin{equation}
\bar{\P} = \diag ( \bar{\Q}_1, \hat{\P} ).
\end{equation}
\end{itemize}
\label{liftRiccati}
\end{algorithm}

Given $\bar{\x}_K$, we can calculate the feedback control by 
(\ref{augmentedU}). Applying this feedback control to 
(\ref{linearTimeInvariant}), we obtain the next state
$\bar{\x}_{K+1}$.  

\begin{remark}
The proposed lift method for discrete linear periodic 
time-varying system can be used for other design methods such as 
$\H_{\infty}$ and robust pole assignment design.
\end{remark}

\section{Implementation and numerical simulation}

In this section, we discuss some details of implementation that
will reduce some computation time comparing to
the directly implementation described in the previous section. 
We also report the test of the proposed algorithm for 
some real world problems discussed in \cite{yang17, yang16}. We 
compare the test result with the ones obtained in \cite{yang17}.

\subsection{Implementation consideration}

The most expensive calculations in Algorithm \ref{liftRiccati} are
the calculation of $\hat{\Q}$, $\hat{\R}$ , and $\hat{\S}$ in
Step 2, and the calculation of $\hat{\R}^{-1}=\tilde{\R}^{-1}$ 
in Step 3. It is easy to check (cf. \cite{gv93}):
(1) direct calculation of $\hat{\Q}$ requires
${\cal O}(2(p-1)^2n^3)+{\cal O}(2(p-1)n^3)+{\cal O}(n^2)$
flops, (2) direct calculation of $\hat{\R}$ requires
${\cal O}(2p(p-1)^2n^2m)+{\cal O}(2p^2(p-1)nm^2)+{\cal O}(p^2m^2)$
flops, (3) direct calculation of $\hat{\S}$ requires
${\cal O}(2(p-1)^2n^3)+{\cal O}(2(p-1)n^3)$ flops,
and (4) directly calculation of $\hat{\R}^{-1}$ requires 
${\cal O}(p^3m^3)$ flops. For extremely large $p$, i.e., very
long period of the system, the majority of the computation 
is the computation of $\hat{\R}$ and $\hat{\R}^{-1}$.

Let $\Q_A= \Q_1^{\frac{1}{2}} \bar{\A}_1 \in \R^{(p-1)n \times n}$ 
and $\Q_B=\Q_1^{\frac{1}{2}} \bar{\B}_1 \in \R^{(p-1)n \times pm}$.
We will use using Matlab notation for sub-matrices.
Since $\bar{\Q}_1$, $\bar{\Q}_2$, and $\bar{\R}$ are positive
diagonal matrices, we can calculate $\hat{\Q}_1$, $\hat{\S}_2$, 
and $\hat{\R}$ in (\ref{let1}) more efficiently as follows:
\begin{itemize}
\item[] for $i=1:(p-1)n$
\item[] $\hspace{0.3in} \Q_A(i,:)= \Q_1^{\frac{1}{2}}(i,i)
 \bar{\A}_1(i,:)$;
\item[] end
\item[] $\hat{\Q} = \Q_A^{\Tr} \Q_A$ 
\item[] for $i=1:n$
\item[] $\hspace{0.3in} \hat{\Q}(i,i) = \hat{\Q}(i,i)
 +\bar{\Q}_2(i,i)$;
\item[] end
\item[] for $i=1:(p-1)n$
\item[] $\hspace{0.3in} \Q_B(i,:)= \Q_1^{\frac{1}{2}}(i,i)
 \bar{\B}_1(i,:)$;
\item[] end
\item[] $\hat{\R} = \Q_B^{\Tr} \Q_B$ 
\item[] for $i=1:pm$
\item[] $\hspace{0.3in} \hat{\R}(i,i) = \hat{\R}(i,i)
 +\bar{\R}(i,i)$;
\item[] end
\item[] $\hat{\S} = \Q_A^{\Tr} \Q_B$ 
\end{itemize}

It is easy to check (cf. \cite{gv93}) the flops for the following 
calculations: (1) the calculation of $\hat{\Q}$ requires 
${\cal O}((p-1)n)+{\cal O}((p-1)n^2)+{\cal O}(2(p-1)^2n^3)
+{\cal O}(n)$ flops, (2) the calculation of $\hat{\R}$ requires
${\cal O}(p(p-1)nm)+{\cal O}(2p^2(p-1)nm^2)+{\cal O}(pm)$,
(3), the calculation of $\hat{\S}$ requires
${\cal O}(2(p-1)pn^2m)$ flops,
(4) this does not save the computation of  $\hat{\R}^{-1}$.

\subsection{Simulation test for the problem of \cite{yang17}}

The first simulation test problem is the spacecraft attitude
control design using only magnetic torques discussed in 
\cite{yang17}. The number of states of this system is $n=6$.
The number of control inputs of this system is $m=3$.
The controllability of this problem is established
in \cite{yang15}. In this simulation test, we use the same 
discrete-time linear periodic model as in \cite{yang17} with the 
same parameters, such as the spacecraft inertia matrix, orbital 
inclination, orbital altitude, weight matrices $\Q$ and $\R$, and 
the same initial conditions.

Using $p=100$, $p=500$, and $p=1000$, we run all three algorithms 
(one proposed in this paper, one proposed in \cite{yang17}, and 
one proposed in \cite{hl94}) for this design and recorded the CPU 
times for all three algorithms. The result is presented in 
Table 1.

\begin{table}[ht]
\begin{center}
    \begin{tabular}{ | c | c | c |c| }
    \hline
    Samples per period & Algorithm 3.1 &  Yang Algorithm \cite{yang17} & Hench-Laub Algorithm \cite{hl94} \\ \hline
    100 & 0.0097 (s) & 0.0757 (s)   & 0.2711 (s)  \\   \hline
    500 & 0.2528 (s) & 1.6042 (s)   & 6.5435 (s)  \\     \hline
   1000 & 4.2821 (s) & 6.3155 (s)   & 25.8996 (s) \\     \hline
    \end{tabular}
\end{center}
\caption{CPU time comparison for problem in \cite{yang17}}
\end{table}

Clearly, the proposed method is significantly cheaper than 
the algorithms in \cite{hl94,yang17}.

\subsection{Simulation test for the problem of \cite{yang16}}

The second simulation test problem is a combined method for the 
spacecraft attitude and desaturation control design using both
reaction wheels and magnetic torques discussed in 
\cite{yang16}. The number of states of this system is $n=9$.
The number of control inputs of this system is $m=6$.
The controllability of problem is guaranteed
because three reaction wheels are assumed to be available.

Using the parameters provided in \cite{yang16}, for $p=100$, 
$p=500$, and $p=1000$, we obtained solutions for the corresponding 
algebraic Riccati equations and recorded the CPU times for 
all three algorithms. The result is presented in Table 2.

\begin{table}[ht]
\begin{center}
    \begin{tabular}{ | c | c | c |c| }
    \hline
    Samples per period & Algorithm 3.1 &  Yang Algorithm \cite{yang17} & Hench-Laub Algorithm \cite{hl94} \\ \hline
    100 & 0.0284 (s) &  0.1120 (s) & 0.3807 (s) \\     \hline
    500 & 3.6376 (s) &  2.5629 (s) & 9.0144 (s) \\     \hline
   1000 & 38.4912 (s) &  10.0629 (s) & 36.0690 (s) \\     \hline
    \end{tabular}
\end{center}
\caption{CPU time comparison for problem in \cite{yang16}}
\end{table}

For this problem, $m=6$ is twice as large as the previous problem,
the proposed algorithm is faster than the algorithms developed
in \cite{hl94,yang17} when the total number of samples in one
period is moderate ($p=100$ samples per period), but when the 
total number of samples in one period increases (to $p=500$ or 
$p=1000$ samples per period), the advantage of the proposed 
algorithm will be lost because the computation of the inverse of 
$\tilde{\R} \in \R^{6000 \times 6000}$ is ${\cal O}(p^3m^3)$
which is very expensive.

\section{Conclusion}

In this paper, we propose a new lifting method to convert the
discrete-time linear periodic system to a discrete-time linear 
time-invariant system. The LQR design for the discrete-time linear
periodic system is then reduced to the LQR design for the
discrete-time linear time-invariant system. By applying the 
standard algorithm to the discrete-time algebraic Riccati 
equation associated with the augmented LTI system and using
the special structure of the augmented LTI system, we derived an 
efficient algorithm for the LQR design for the discrete time 
linear periodic system. We show that the new algorithm is very
efficient compared to the existing ones 
for the discrete-time periodic algebraic Riccati equation, in 
particular when the number of samples in one period is moderate. 
We demonstrated this property by the numerical test for the
design problems of spacecraft attitude control using magnetic
torques. The proposed lift method for discrete linear 
periodic time-varying system can be used for other design
methods such as $\H_{\infty}$ and robust pole assignment designs.

%

\end{document}